\documentclass{amsart}

\usepackage{color}
\usepackage{epsfig}
\usepackage{amsthm}
\usepackage{amssymb}
\usepackage{amsmath}
\usepackage{amscd}

\newtheorem {Theorem}  {Theorem}

\numberwithin{Theorem}{section}

\newtheorem {Lemma}[Theorem]  {Lemma}
\newtheorem {Proposition}[Theorem]{Proposition}
\theoremstyle{definition}
\newtheorem{Definition}[Theorem]{Definition}
\theoremstyle{remark}
\newtheorem{Remark}[Theorem]{Remark}
\newtheorem{Example}[Theorem]{Example}

\newtheorem {Corollary}[Theorem]{Corollary}

\expandafter\chardef\csname pre amssym.def
at\endcsname=\the\catcode`\@ \catcode`\@=11
\def\undefine#1{\let#1\undefined}
\def\newsymbol#1#2#3#4#5{\let\next@\relax
 \ifnum#2=\@ne\let\next@\msafam@\else
 \ifnum#2=\tw@\let\next@\msbfam@\fi\fi
 \mathchardef#1="#3\next@#4#5}
\def\mathhexbox@#1#2#3{\relax
 \ifmmode\mathpalette{}{\m@th\mathchar"#1#2#3}%
 \else\leavevmode\hbox{$\m@th\mathchar"#1#2#3$}\fi}
\def\hexnumber@#1{\ifcase#1 0\or 1\or 2\or 3\or 4\or 5\or 6\or 7\or 8\or
 9\or A\or B\or C\or D\or E\or F\fi}

\font\teneufm=eufm10 \font\seveneufm=eufm7 \font\fiveeufm=eufm5
\newfam\eufmfam
\textfont\eufmfam=\teneufm \scriptfont\eufmfam=\seveneufm
\scriptscriptfont\eufmfam=\fiveeufm

\catcode`\@=\csname pre amssym.def at\endcsname

\newcounter{remark}
\def  \eps  {\epsilon}

\def\3n{\negthinspace \negthinspace \negthinspace }
\def\2n{\negthinspace \negthinspace }
\def\1n{\negthinspace }

\newcommand{\bg}{\begin{equation}}
\newcommand{\ed}{\end{equation}}
\newcommand{\bga}{\begin{eqnarray}}
\newcommand{\eda}{\end{eqnarray}}

\def\cbdu{\hfill{$\Box$}}

\renewcommand{\a}{\alpha}

\newcommand{\R}{\mathbf{R}}

\def\bex{\begin{Example}}
\def\eex{\end{Example}}

\def  \R   {{\mathbb R}}

\def  \12  {{\frac{1}{2}}}

\def\bd{\begin{Definition}}
\def\ede{\end{Definition}}
\def\be{\begin{equation}}
\def\bel{\begin{equation}\label}
\def\ee{\end{equation}}
\def\bt{\begin{Theorem}}
\def\et{\end{Theorem}}
\def\bc{\begin{Corollary}}
\def\ec{\end{Corollary}}
\def\bl{\begin{Lemma}}
\def\el{\end{Lemma}}
\def\bp{\begin{Proposition}}
\def\ep{\end{Proposition}}

\def\br{\begin{Remark}}
\def\er{\end{Remark}}
\def\ba{\begin{array}}
\def\ea{\end{array}}

\def\bea{\begin{eqnarray}}
\def\eea{\end{eqnarray}}

\begin{document}

\title[Decay of solutions to Navier-Stokes-Voigt equations]{Decay characterization of solutions to Navier-Stokes-Voigt equations in terms of the initial datum}

\author{C\'esar J. Niche}
\address[C.J. Niche]{Departamento de Matem\'atica Aplicada, Instituto de Matem\'atica, Universidade Federal do Rio de Janeiro, CEP 21941-909, Rio de Janeiro - RJ, Brasil}
\email{cniche@im.ufrj.br}





\begin{abstract}
The Navier-Stokes-Voigt equations are a regularization of the Na\-vier-Stokes equations that share some of its asymptotic and statistical properties and have been used in direct numerical simulations of the latter. In this article we characterize the decay rate of solutions to the Navier-Stokes-Voigt equations in terms of the decay character of the initial datum and study the long time behaviour of its solutions by comparing them to solutions to the linear part. 
\end{abstract}


\maketitle

\section{Introduction}

The Navier-Stokes  equations in $\R^3$

\begin{eqnarray}
\label{eqn:navier-stokes}
\partial_t u - \nu \Delta u  +  (u \cdot \nabla) u +  \nabla p & = & 0, \nonumber \\ \nabla \cdot u & = & 0, \nonumber \\ u_0 (x) & = & u(x,0)
\end{eqnarray}
with $\nu > 0$, describe the time evolution of the velocity $u$ and the pressure $p$ of a homogeneous incompressible viscous fluid and are of great  importance in the study of turbulent flows. The analytical and numerical study of these equations is very difficult and as a result of this, many different kinds of approximations and regularizations  have been proposed in order to have a better understanding of  the behaviour of its solutions. 

One such regularization that has been extensively studied in recent years, in the case of bounded or periodic domains, is given by the Navier-Stokes-Voigt  equations

\begin{eqnarray}
\label{eqn:navier-stokes-voigt}
\partial_t \left(u - \a ^2 \Delta  u \right) - \nu \Delta u  +  (u \cdot \nabla) u +  \nabla p & = & 0, \nonumber \\ \nabla \cdot u & = & 0, \nonumber \\ u_0 (x) & = & u(x,0),
\end{eqnarray}
supplemented with appropiate boundary conditions. These equations were originally used by Oskolkov \cite{MR0377311} to study the motion of viscoelastic flows with relaxation time, i.e. the time that the flow takes to respond to an external force, of order $\a^2 / \nu$ and have been proposed by Cao, Lunasin and Titi \cite{MR2264822} as a good model for numerical simulations of the Navier-Stokes equations, provided $\a > 0$ is small.

The introduction of the regularizing term $- \a ^2 \partial_t \Delta  u$ in (\ref{eqn:navier-stokes-voigt}) has many important consequences. We first note that, assuming $u$ is a regular enough solution, (\ref{eqn:navier-stokes-voigt}) has a natural energy equality

\be
\label{eqn:energy-equality-nsv}
\frac{1}{2} \frac{d}{dt} \left(  \Vert u(t) \Vert ^2 _{L^2 (\R^3)} + \a^2 \Vert \nabla u(t) \Vert ^2 _{L^2 (\R^3)} \right) = - \nu \Vert \nabla u(t) \Vert ^2 _{L^2 (\R^3)},
\ee
instead of the usual one for (\ref{eqn:navier-stokes}). Moreover, Navier-Stokes-Voigt equations do not have a parabolic character, as Navier-Stokes equations do, behaving instead as a damped hyperbolic system. This leads to wellposedness for equation (\ref{eqn:navier-stokes-voigt}) both forward and backwards in time (see Cao, Lunasin and Titi \cite{MR2264822}, Kalantarov and Titi \cite{MR2570790}) and to differences in the energy spectrum when compared to (\ref{eqn:navier-stokes}) (see Levant, Ramos and Titi \cite{MR2655910}, Ramos and Titi \cite{MR2629486}). However, solutions to Navier-Stokes-Voigt equations share many other statistical and asymptotic properties with the Navier-Stokes equations when $\a > 0$ is small (see  Kalantarov, Levant and Titi \cite{MR2495891}, Larios and Titi \cite{Larios2010603}, Levant, Ramos and Titi \cite{MR2655910}, Ramos and Titi \cite{MR2629486}), while presenting advantages for the implementation of numerical models for their study. This supports the idea of using (\ref{eqn:navier-stokes-voigt}) in direct numerical simulations of (\ref{eqn:navier-stokes}) (see Borges and Ramos \cite{MR3048413}, Kuberry, Larios, Rebholz and Wilson \cite{MR2970841}, Layton and Rebholz \cite{MR3171807}). For more results concerning the Navier-Stokes-Voigt equations, see Bersellli and Bisconti \cite{MR2846786}, Ebrahimi, Holst and Lunasin \cite{MR3116155}, Garc\'{\i}a-Luengo, Mar\'{\i}n-Rubio and Real \cite{MR2891552}, Gao and Sun \cite{MR2863948}, Tang \cite{MR3217168}, Yue and Zhong {\cite{MR2806333} and references therein. 

Only recently, Zhao and Zhu \cite{MR3316059}  adressed questions concerning the Navier-Stokes-Voigt equations in the whole space $\R^3$. In their work, they proved  existence of weak solutions  to  (\ref{eqn:navier-stokes-voigt}), but their main result  concerns the decay rate for the $H^1$ norm of solutions, which, for  initial data $u_0$ in $L^1 (\R^3) \cap H^1 (\R^3)$, is proved to be

\be
\label{eqn:decay-zhao-zhu}
\Vert u(t) \Vert _{H^1 (\R^3)} ^2 \leq C (1 + t) ^{- \frac{3}{2}}.
\ee
The main tool in the proof of these estimates is a clever adaptation of the Fourier Splitting method of M.E. Schonbek \cite{MR571048}, \cite{MR775190}, \cite{MR837929}, which takes into account the extra linear term present in (\ref{eqn:navier-stokes-voigt}) with respect to (\ref{eqn:navier-stokes}). 

The main goal of this article is to prove a decay estimate for the Navier-Stokes-Voigt equations for {\em any} initial datum in $H^1 _{\a} (\R^3)$, where

\begin{displaymath}
\Vert u \Vert _{H^1 _{\a} (\R^3)} ^2 = \Vert u \Vert _{L^2 (\R^3)} ^2 + \a^2 \Vert \nabla u \Vert _{L^2 (\R^3)} ^2 
\end{displaymath}
is the natural norm associated to (\ref{eqn:energy-equality-nsv}). To achieve this, we associate to every $u_0$ in $L^2 (\R^3)$ a decay character $r^{\ast} = r^{\ast} (u_0)$, which measures the ``order'' of $\widehat{u_0} (\xi)$ at $\xi = 0$ in frequency space and then find sharp decay estimates for solutions to the linear part

\be
\label{eqn:linear-part}
\partial_t \left( u - \a^2 \Delta u \right) - \nu \Delta u  = 0
\ee
in terms of $r^{\ast}$. We then use the Fourier Splitting Method and these estimates for the linear part to characterize the decay of Navier-Stokes-Voigt equations. This same strategy was followed by Niche and M.E. Schonbek \cite{2015arXiv150102105N}, \cite{Niche01042015} to obtain analogous results for the Navier-Stokes equations and some parabolic compressible approximations and for the dissipative quasi-geostrophic equation.

We state now one of the main results of this article.

\bt
\label{main-theorem}
Let $u_0 \in H^1 _{\a} (\R^3)$, with $\nabla \cdot u_0 = 0$ and $r^{\ast} = r^{\ast} (u_0)$, $- \frac{3}{2} \leq r^{\ast} < \infty$. Then for $\nu > \a^2$, every weak solution to (\ref{eqn:navier-stokes-voigt}) is such that for $t > 0$

\begin{displaymath}
\Vert u(t) \Vert _{H^1 _{\a}} ^2 \leq C_1 (C_2 + t) ^{- \min \{\frac{3}{2} + r^{\ast}, \frac{5}{2} \}},
\end{displaymath}
for some positive constants $C_i = C_i (\Vert u_0 \Vert _{H^1 _{\a} (\R^3)}, \a, \nu)$, $i = 1, 2$.
\et

The decay in Theorem \ref{main-theorem} is similar to the one for the Navier-Stokes equations

\be
\label{eqn:characterization-ns}
\Vert u(t) \Vert _{L^2 (\R^3)} ^2 \leq C_1 (C_2 + t) ^{- \min \{ \frac{3}{2} + r^{\ast}, \frac{5}{2} \}}, \qquad - \frac{3}{2} \leq r^{\ast} < \infty,
\ee
obtained by Bjorland and M.E. Schonbek \cite{MR2493562} (see also Niche and M.E. Schonbek \cite{Niche01042015}).  In our work, using the Fourier Splitting Method we obtain in frequency space the differential inequality (\ref{eqn:previous-to-main-inequality}), in which the decay rate is bounded from above by the average of the solution in a shrinking, time-dependent ball containing the small frequencies. The leading terms in the right hand side of this expression are the same as in the Navier-Stokes equations and, as the other terms provide faster decay rates, we recover the same estimate.

The bound from Theorem \ref{main-theorem} allows us to determine the decay rate for the difference between the whole solution $u(t)$ and the solution $\overline{u} (t)$ to the linear part (\ref{eqn:linear-part}) with the same initial datum. 

\bt
\label{decay-nonlinear-part}
Let $u_0 \in H^1 _{\a} (\R^3)$, with $\nabla \cdot u_0 = 0$ and $r^{\ast} = r^{\ast} (u_0)$, $- \frac{3}{2} \leq r^{\ast} < \infty$. Let $w = u - \overline{u}$. Then for $\nu > \a^2$,

\begin{displaymath}
\Vert w(t) \Vert ^2 _{H^1 _{\a} (\R^3)} \leq C_1 (C_2 + t) ^{- \min \{ \frac{9}{4} + \frac{3}{2} r^{\ast}, \frac{13}{4} + \frac{1}{2} r^{\ast}, \frac{11}{2}  \}}.
\end{displaymath}
for some positive constants $C_i = C_i (\Vert u_0 \Vert _{H^1 _{\a} (\R^3)}, \a, \nu)$, $i = 1, 2$.
\et

This result has quantitative and qualitative differences with the analogous one for the $L^2 (\R^3)$ norm of $w = u - \overline{u}$ in the Navier-Stokes equations, where

\begin{displaymath}
\Vert w(t) \Vert _{L^2 (\R^3)} ^2 \leq C_1 (C_2 + t) ^{- \min \{\frac{7}{4} + r^{\ast}, \frac{5}{2} \}},
\end{displaymath}
see Niche and M.E. Schonbek \cite{Niche01042015} \footnote{We note that in \cite{Niche01042015} this estimate was proved for a compressible approximation to the Navier-Stokes equations. However, the proof carries over with no modifications due to the very similar structure of these equations, regardless of compressibility. For a similar result for Navier-Stokes, see Wiegner \cite{MR881519}.} .The decay rate from Theorem \ref{decay-nonlinear-part} is faster for $r^{\ast} > -1$, but for $r^{\ast} = -\frac{3}{2}$, the difference $w$ in the Navier-Stokes-Voigt equations may have arbitrarily slow decay, while in the Navier-Stokes equations we always have  some uniform decay with speed at least $(1 + t) ^{- \frac{1}{4}}$. We address this interesting point in Remark \ref{remark-slow-decay-nsv}.

This article is organized as follows. In Section \ref{section-decay-character} we recall the definition of the decay character of a function in $H^s (\R^n)$, with $s \geq 0$ and we prove Theorem \ref{characterization-decay-linear-part}, in which we caracterize the decay of solutions to the linear system (\ref{eqn:linear-part}), for initial data in $H^1 _{\a} (\R^n)$. Then, in Section \ref{section-proofs}, we prove Theorems \ref{main-theorem} and \ref{decay-nonlinear-part}. 

\subsection*{Acknowledgments}
The author thanks Fabio Ramos and Mar\'{\i}a E. Schonbek for their comments and remarks on parts of this article.

\section{Decay character and characterization of decay of linear part}
\label{section-decay-character}

\subsection{Definition and properties of decay character}

The Fourier Splitting Method, introduced by M.E. Schonbek \cite{MR571048}, \cite{MR775190}, \cite{MR837929} to study decay of viscous conservation laws and the Navier-Stokes equations, relies on the study of small frequencies of the solution  to estimate the decay rates. In order to prove sharp estimates for the decay of the heat equation in terms of the initial data, Bjorland and M.E. Schonbek \cite{MR2493562} introduced the idea of {\em decay character $r^{\ast} = r^{\ast} (u_0)$} of a function $u_0$ in $L^2 (\R^n)$.  Roughly speaking, the decay character is the ``order'' of $| \widehat{u_0} (\xi)|$ in frequency space at $\xi = 0$. With these estimates in hand, Bjorland and M.E. Schonbek characterized the decay of solutions to the Navier-Stokes equations in terms of $r^{\ast}$, see (\ref{eqn:characterization-ns}).

Recently, Niche and M.E. Schonbek \cite{Niche01042015} generalized this notion in order to use data in $H^s (\R^n)$, for $s  > 0$ and to obtain results for other linear systems. We recall now some of these definitions and results.

\bd Let  $u_0 \in L^2(\R^n)$ and $\Lambda = (- \Delta) ^{\frac{1}{2}}$. For $s \geq 0$ the {\em s-decay indicator}  $P^s _r (u_0)$ corresponding to $\Lambda ^s u_0$ is 

\begin{displaymath}
P^s _r (u_0) = \lim _{\rho \to 0} \rho ^{-2r-n} \int _{B(\rho)} |\xi|^{2s}| \widehat{u}_0 (\xi)| ^2 \, d \xi
\end{displaymath}
for $r \in \left(- \frac{n}{2} + s, \infty \right)$, where $B(\rho)$ is the ball at the origin with radius $\rho$.
\ede

\br \label{heuristics} \, Setting $r = q + s$, we see that the $s$-decay indicator compares $|\widehat{\Lambda ^s u_0} (\xi)|^2$ to $f(\xi) = |\xi|^{2(q + s)}$ near $\xi = 0$. When $s = 0$ we recover the definition of Bjorland and M.E. Schonbek \cite{MR2493562}.
\er

\bd  \label{df-decay-character} The {\em decay character of $\Lambda^s u_0$}, denoted by $r_s ^{\ast} = r_s ^{\ast}( u_0)$ is the unique  $r \in \left( -\frac{n}{2} + s, \infty \right)$ such that $0 < P^s _r (u_0) < \infty$, provided that this number exists. If such  $P^s _r ( u_0)$ does not exist, we set $r_s ^{\ast} = - \frac{n}{2} + s$, when $P^s _r (u_0)  = \infty$ for all $r \in \left( - \frac{n}{2} + s, \infty \right)$  or $r_s ^{\ast} = \infty$, if $P^s _r (u_0)  = 0$ for all $r \in \left( -\frac{n}{2} + s, \infty \right)$. 
\ede

\br \label{remark-lp-l2} \, Let $u_0 \in L^2 (\R^n)$ such that $\widehat{u}_0 (\xi) = 0$, for $|\xi| < \delta$, for some $\delta > 0$. Then, $P_r ^s (u_0) = 0$, for any $r \in \left( - \frac{n}{2} + s, \infty\right)$ and  $r^{\ast} _s (u_0) = \infty$.  
\er

\br If $u_0 \in L^p (\R^n) \cap L^2 (\R^n)$, with $1  \leq p \leq 2$, then $r^{\ast} (u_0) = - n \left( 1 - \frac{1}{p} \right)$, so if $u_0 \in L^1 (\R^n) \cap L^2 (\R^n)$ we have that $r^*(u_0) =0$ and if $u_0 \in L^2 (\R^n)$ but $u_0 \notin L^p (\R^n)$, for any $1 \leq p < 2$, we have that $r^*(u_0) = -  \frac{n}{2}$. For details, see Niche and M.E. Schonbek \cite{Niche01042015}.
\er

\smallskip

Let $u_0 \in H^s (\R^n)$, with $s > 0$. As $\widehat{\Lambda^s (u_0)} (\xi) = |\xi|^s \widehat{u_0} (\xi)$, the heuristics given in Remark \ref{heuristics} lead us to expect that $r^{\ast} _s (u_0) = r^{\ast} (\Lambda ^s u_0) = s + r^{\ast} (u_0)$. This is the content of the following Theorem.

\bt \label{decay-character-hs} (Theorem 2.11, Niche and M.E. Schonbek \cite{Niche01042015}) Let $u_0 \in H^s (\R^n)$, with $s > 0$.  
\begin{enumerate}
 \item If $-\frac{n}{2} < r^{\ast} (u_0) < \infty$ then  $- \frac{n}{2} +s< r_s^{\ast}(u_0) < \infty$ and  $r_s^{\ast}(u_0) = s + r^{\ast} (u_0)$. 
\item $r_s^{\ast}  (u_0) = \infty$ if and only if $r^{\ast} (u_0) = \infty$.
\item  $r^{\ast} (u_0) =- \frac{n}{2}$ if and only if $r_s^{\ast}(u_0) = r^{\ast} (u_0) + s = - \frac{n}{2} + s$. \end{enumerate}
\et

\subsection{Characterization of decay of linear part}  We address now the system

\be
\label{eqn:linear-part-nsv}
\partial_t \left(v - \a ^2 \Delta  v \right) - \nu \Delta v= 0,
\ee
which is the linear part of the Navier-Stokes-Voigt equations (\ref{eqn:navier-stokes-voigt}). From it, we obtain the energy equality 

\be
\label{eqn:energy-inequality-linear-part}
\frac{d}{dt} \left( \Vert v(t) \Vert ^2 _{L^2 (\R^n)} + \a ^2 \Vert \nabla v(t) \Vert ^2 _{L^2 (\R^n)} \right) + 2 \nu \Vert \nabla v(t) \Vert _{L^2 (\R^n)} ^2 = 0
\ee
which naturally leads us to study the behaviour of the $H^1 _{\a}$ norm

\begin{displaymath}
\Vert v \Vert _{H^1 _{\a} (\R^n)} ^2 = \Vert v \Vert _{L^2 (\R^n)} ^2 + \a^2 \Vert \nabla v \Vert _{L^2 (\R^n)} ^2. 
\end{displaymath}
Then in Fourier space, the solution to (\ref{eqn:linear-part-nsv}) with initial datum $v_0$ is 

\begin{displaymath}
\widehat{v} (\xi, t) = e^{t \mathcal{M} (\xi)} \widehat{v_0} (\xi)
\end{displaymath}
where
\be
\label{eqn:nsv-multiplier}
e^{t \mathcal{M} (\xi)} = e^{ - \frac{\nu |\xi|^2}{1 + \a^2 |\xi|^2} t}.
\ee

\br
\label{facts-linear-part-nsv}
Note that (\ref{eqn:nsv-multiplier}) is bounded but not in any $L^p$, for $1 \leq p < \infty$. Then, solutions to (\ref{eqn:linear-part-nsv}) with initial data in Sobolev space remain there but are not regularized, i.e. they do not belong to higher order Sobolev spaces, as the kernel associated to (\ref{eqn:nsv-multiplier}) does not ``absorb'' derivatives as the heat kernel does. However, for small frequencies (\ref{eqn:nsv-multiplier}) behaves like a Gaussian, which allows us to the Fourier Splitting Method.
\er

In the following Theorem, we characterize the decay of solutions to (\ref{eqn:linear-part-nsv}) in terms of the decay character $r^{\ast} = r^{\ast} (v_0)$.

\bt
\label{characterization-decay-linear-part}
Let $v_0 \in H^1 _{\alpha} (\R^n)$, with $r^{\ast} (v_0) = r^{\ast}$. If $\nu > \a^2$, then for the solution $v(t)$ to (\ref{eqn:linear-part-nsv}) with initial datum $v_0$ we have

\begin{enumerate}
\item if $- \frac{n}{2} < r^{\ast} < \infty$, then there exist constants $C_1, C_2, C_3 > 0$ such that

\begin{displaymath}
C_1 (1 + t) ^{- \left( \frac{n}{2} + r^{\ast} \right)} \leq \Vert v(t) \Vert ^2 _{H^1 _{\alpha} (\R^n)} \leq C_2 (C_3 + t) ^{- \left( \frac{n}{2} + r^{\ast} \right)};
\end{displaymath}
\item if $r^{\ast} = - \frac{n}{2}$ given any $\eps > 0$, there exists $C = C(\eps) > 0$ such that

\begin{displaymath}
\Vert v(t) \Vert ^2 _{H^1 _{\alpha} (\R^n)} \geq C (1 + t) ^{- \eps}, \quad \forall \, t > 0,
\end{displaymath}
i.e. the decay is slower than any algebraic rate;
\item if $r^{\ast} = \infty$, given any $m > 0$ then there exists $C = C(m) > 0$ such that

\begin{displaymath}
\Vert v(t) \Vert ^2 _{H^1 _{\alpha} (\R^n)} \leq C (1 + t) ^{-m}, \quad \forall \, t  > 0,
\end{displaymath}
i.e. the decay is faster than any algebraic rate.
\end{enumerate}
\et

\medskip

{\bf Proof:} (1) Recall first that we have, by Theorem \ref{decay-character-hs}, that $r^{\ast} (\nabla v_0) = 1 + r^{\ast} (v_0)$. As $P_r (v_0) > 0$ there exist $\rho_0 > 0$ and $C_1 > 0$ such that for $0 < \rho < \rho_0$ 

\be
\label{eqn:decay-indicator-below}
C_1 < \rho ^{-2r - n} \int _{B(\rho)} |\widehat{v_0} (\xi)|^2 \, d \xi, \qquad C_1 < \rho ^{-2(r + 1) - n} \int _{B(\rho)} |\xi|^2 |\widehat{v_0} (\xi)|^2 \, d \xi.
\ee
Consider $B(t) = \{ \xi \in \R^n: |\xi| \leq \rho(t) \}$, for some continuous, positive and decreasing $\rho$ to be determined later. Then 

\begin{eqnarray}
\label{eqn:inequality-lower-bound}
\Vert v(t) \Vert _{H^1 _{\alpha}} ^2 & \geq & \int _{B(t)} \left( 1 + \alpha ^2 |\xi|^2 \right) |e^{t \mathcal{M} (\xi)} \widehat{v_0} (\xi)|^2 \, d \xi \nonumber \\ & \geq & \int _{B(t)} \left( 1 + \alpha^2  |\xi|^2 \right) e^{- 2\nu |\xi|^2 t} |\widehat{v_0} (\xi)|^2 \, d \xi \nonumber \\ & \geq & C \rho ^{2r + n} (t) e^{- 2\nu \rho ^2 (t)} \rho ^{- 2r - n} (t) \int _{B(t)} |\widehat{v_0} (\xi)|^2 \, d \xi \nonumber \\ & + & C \rho ^{2(r + 1) + n} (t) e^{- 2\nu \rho ^2 (t)} \rho ^{- 2(r + 1) - n} (t) \int _{B(t)} |\xi|^2 |\widehat{v_0} (\xi)|^2 \, d \xi \nonumber \\ & \geq & C \rho ^{2r + n} (t) e^{- 2\nu \rho ^2 (t)} + C \rho ^{2(r + 1) + n} (t) e^{- 2\nu \rho ^2 (t)}
\end{eqnarray}
where we used (\ref{eqn:decay-indicator-below}) in the last line. Taking $\rho (t) = \rho_0 (1 + t) ^{-\frac{1}{2}}$ we obtain

\begin{displaymath}
\Vert v(t) \Vert _{H^1 _{\alpha}} ^2  \geq  C (1 + t) ^{- \left( r + \frac{n}{2} \right)} + C (1 + t) ^{- \left( r + 1 + \frac{n}{2} \right)} \geq C (1 + t) ^{- \left( r + \frac{n}{2} \right)},
\end{displaymath}
from which the result follows. To obtain the upper bound we use the Fourier Splitting method. We  follow the ideas in Zhu and Zhao \cite{MR3316059}. From the energy inequality (\ref{eqn:energy-inequality-linear-part})  we obtain

\begin{displaymath}
\frac{d}{dt} \left( \int _{\R^n} \left( 1 + \a^2 |\xi|^2 \right) |\widehat{u} (\xi, t) |^2 \, d \xi \right) + 2 \nu \int _{\R^n}  |\xi|^2  |\widehat{u} (\xi, t) |^2 \, d \xi = 0.
\end{displaymath}
Let 

\begin{displaymath}
B(t) = \left\{\xi \in \R^n: |\xi| \leq \rho (t) = \left( \frac{g'(t)}{2 \nu g(t) - \a^2 g'(t)} \right) ^{\frac{1}{2}} \right\}
\end{displaymath}
where $g \in C^1 (\R)$, $g(0) = 1$, $g'(t) > 0$ and $2 \nu g(t) - \a^2 g'(t) > 0$, for all $t > 0$. We then have

\begin{align}
\label{eqn:inequality-before-decay-indicator}
\frac{d}{dt}  \left( g(t) \int _{\R^n} \left( 1 + \a^2 |\xi|^2 \right)  |\widehat{v} (\xi, t) |^2 \, d \xi \right)  \leq   g'(t) \int _{B(t)} \left( 1 + \a^2 |\xi|^2 \right) |\widehat{v} (\xi, t) |^2 \, d \xi \nonumber \\  \leq  g'(t) \rho^{2r + n} (t) \rho^{- 2r - n} (t) \int _{B(t)} |\widehat{v_0} (\xi)|^2 \, d \xi \nonumber \\  +   \a^2 g'(t) \rho^{2(r + 1) + n} (t) \rho^{- 2(r + 1) - n} (t) \int _{B(t)} |\xi|^2 |\widehat{v_0} (\xi)|^2 \, d \xi
\end{align} 
where the first inequality was proved in page 186 in Zhu and Zhao \cite{MR3316059}. As $P_r (u_0) < \infty$ and $r^{\ast} (\nabla u_0) = 1 + r^{\ast} (u_0)$, there exist $\rho_0 > 0$ and $C > 0$ such that for $0 < \rho < \rho _0$

\be
\label{eqn:decay-indicator-above}
\rho ^{-2r - n} \int _{B(\rho)} |\widehat{v_0} (\xi)|^2 \, d \xi \leq C, \qquad \rho ^{-2(r + 1) - n} \int _{B(\rho)} |\xi|^2 |\widehat{v_0} (\xi)|^2 \, d \xi \leq C.
\ee
Then, from (\ref{eqn:inequality-before-decay-indicator}) and (\ref{eqn:decay-indicator-above}) we obtain 

\be
\label{eqn:previous-to-estimate}
\frac{d}{dt}  \left( g(t) \int _{\R^n} \left( 1 + \a^2 |\xi|^2 \right) |\widehat{u} (\xi, t) |^2 \, d \xi \right) \leq C g'(t) \rho^{2r + n} (t) + C g'(t) \rho^{2(r + 1) + n} (t).
\ee
Taking $g(t) = \left( \frac{1}{2 \nu} t + b \right)^k$, where $k > 1 + r + \frac{n}{2}$ and $b > \frac{k \a^2}{2 \nu ^2}$, we have that $2 \nu g(t) - \a^2 g'(t) > 0$ holds (because $\nu > \a^2$) and that $\rho(t) = C (t + b) ^{- \frac{1}{2}}$. From (\ref{eqn:previous-to-estimate}) the estimate 

\begin{displaymath}
\Vert v(t) \Vert _{H^1 _{\alpha}} ^2 \leq C (t + b) ^{-k} + C (t + b) ^{- \left( r + \frac{n}{2} \right)} + C (t + b) ^{- \left( r + 1 + \frac{n}{2} \right) } \leq  C (t + b) ^{- \left( r + \frac{n}{2} \right)}
\end{displaymath}
follows.

\medskip
\noindent (2) As $r^{\ast} = - \frac{n}{2}$, for any  $r \in (- \frac{n}{2}, \infty)$ we have that $P_r (v_0) = \infty$. Then for any constant $\widetilde{C} = \widetilde{C} (r)$ there exists $\rho_0 > 0$ such that for $0 < \rho_0 < \rho$

\begin{displaymath}
\widetilde{C} < \rho ^{-2r - n} \int _{B(\rho)} |\widehat{v_0} (\xi)|^2 \, d \xi, \qquad \widetilde{C} < \rho ^{-2 (r + 1) - n} \int _{B(\rho)} |\xi|^2 |\widehat{v_0} (\xi)|^2 \, d \xi.
\end{displaymath}
As in the proof of the lower bound in (1) with an inequality similar to (\ref{eqn:inequality-lower-bound}), we obtain

\begin{displaymath}
C (1 + t)^{-r - \frac{n}{2}} \leq \Vert v(t) \Vert _{H^1 _{\alpha}} ^2.
\end{displaymath}
As this holds for any $r \in (- \frac{n}{2}, \infty)$, the estimate is established.

\noindent (3) Since $r^{\ast} = \infty$, for any $r \in (- \frac{n}{2}, \infty)$ we have $P_r (v_0) = 0$. So, for any $\widetilde{C} = \widetilde{C} (r)$ there exists $\rho_0 > 0$ such that for $0 < \rho_0 < \rho$

\begin{displaymath}
\rho ^{-2r - n} \int _{B(\rho)} |\widehat{v_0} (\xi)|^2 \, d \xi < \widetilde{C}, \rho ^{-2(r + 1) - n} \int _{B(\rho)} |\xi|^2 |\widehat{v_0} (\xi)|^2 \, d \xi < \widetilde{C}
\end{displaymath}
Proceeding as in the proof of the upper bound in (1)   with an inequality similar to (\ref{eqn:previous-to-estimate}), we obtain

\begin{displaymath}
\Vert v(t) \Vert _{H^1 _{\alpha}} ^2 \leq C (1 + t)^{-r - \frac{n}{2}} .
\end{displaymath}
As this holds for any $r \in (- \frac{n}{2}, \infty)$, the estimate is established. $\Box$

\section{Proof of Theorems \ref{main-theorem} and \ref{decay-nonlinear-part}}
\label{section-proofs}

Before proceeding to the proof of  Theorem \ref{main-theorem} we prove some Lemmas. Suppose that the solution to the Navier-Stokes-Voigt equation is regular enough such that

\be
\label{eqn:solution-fourier-space}
\widehat{u} (\xi, t) = \widehat{u_0} (\xi) \, e^{t \mathcal{M} (\xi)} - \int _0 ^t e^{(t -\tau) \mathcal{M} (\xi)} G(\xi, \tau) \, d \tau,
\ee
where $G(\xi,t) =  \mathcal{F} \left( (u \cdot \nabla) u - \nabla p \right) (\xi, t)$ and $\mathcal{F}$ denotes the Fourier transform. 

\begin{Lemma} \label{estimate-modulus-uhat} For $\widehat{u}$ as in (\ref{eqn:solution-fourier-space}) we have

\be
\label{eqn:estimate-modulus-solution-fourier}
|\widehat{u} (\xi, t)|^2 \leq C  e^{2t \mathcal{M} (\xi)} |\widehat{u_0} (\xi)|^2  + C |\xi|^2 \left( \int _0 ^t \Vert u(\tau) \Vert _{L^2 (\R^3)} ^2  \, d \tau \right)^2.
\ee
\end{Lemma}

{\bf Proof:} The estimate follows from 

\begin{eqnarray*}
|\mathcal{F} ((u \cdot \nabla) u) (\xi,t)| & \leq & |\mathcal{F} (\nabla \cdot (u \otimes u)) (\xi, t)| \leq C |\xi| |\mathcal{F} (u \otimes u)) (\xi, t)| \\ & \leq & C |\xi| \Vert \mathcal{F} (u \otimes u) (t) \Vert _{L^{\infty} (\R^3)} \nonumber \\ & \leq & C |\xi| \Vert u \otimes u (t) \Vert _{L^1 (\R^3)} \leq C |\xi| \Vert u(t) \Vert _{L^2 (\R^3)} ^2
\end{eqnarray*}
and from $\Delta p = - \nabla (u \cdot \nabla) u$, as 

\begin{displaymath}
|\mathcal{F} (\nabla p) (\xi, t)| \leq |\mathcal{F} ((u \cdot \nabla) u) (\xi, t)| \leq C |\xi| \Vert u(t) \Vert _{L^2 (\R^3)} ^2. \qquad \Box
\end{displaymath}

{\bf Proof of Theorem \ref{main-theorem}:} \,  We again use the Fourier Splitting Method. We recall that existence of weak solutions to (\ref{eqn:navier-stokes-voigt}) has been established by Zhu and Zhao, see Theorem 2.1 in  \cite{MR3316059}. As usual in this method, we prove the estimates for regular enough solutions and we then pass to the limit  to establish it for weak solutions. Full details for this argument can be found in pages 267-269 in Lemarie-Rieusset \cite{MR1938147} and in the Appendix in Wiegner \cite{MR881519}. 

 We follow the ideas of Zhu and Zhao \cite{MR3316059}.  As before, let 

\be
\label{eqn:ball}
B(t) = \left\{ \xi \in \R^3: |\xi| \leq \rho (t) = \left( \frac{g'(t)}{2 \nu g(t) - \a^2 g'(t)} \right) ^{\frac{1}{2}} \right\}
\ee
where $g \in C^1 (\R)$, $g(0) = 1$, $g'(t) > 0$ and $2 \nu g(t) - \a^2 g'(t) > 0$, for all $t > 0$. The incompressibility condition leads to the energy inequality (\ref{eqn:energy-equality-nsv}), which is the same one as for the linear part, see (\ref{eqn:energy-inequality-linear-part}). We then obtain the same inequality as in the first line in (\ref{eqn:inequality-before-decay-indicator}), which after integration and use of Lemma \ref{estimate-modulus-uhat}  becomes

\begin{align}
\label{eqn:previous-to-main-inequality}
 g(t) \int _{\R^3} \left( 1 + \a^2 |\xi|^2 \right) |\widehat{u} (\xi, t) |^2 \, d \xi   \leq   C + \int _0 ^t g'(s) \int _{B(s)} \left( 1 + \a^2 |\xi|^2 \right) |\widehat{u} (\xi, s) |^2 \, d \xi \nonumber \\ \leq   C + C \int _0 ^t g'(s) \left( \int _{B(s)} \left( 1 + \a^2 |\xi|^2 \right) e^{2 t \mathcal{M} (\xi)} |\widehat{u_0} (\xi)|^2  \, d \xi \right) \, ds \nonumber \\  +   C \int _0 ^t  g'(s)  \left( \int _{B(s)} \left(1  +  \a^2 |\xi|^2 \right)  |\xi|^2 \left( \int _0 ^t \Vert u(\tau) \Vert _{L^2 (\R^3)} ^2  \, d \tau \right)^2  \, d \xi \right) \, ds.
\end{align}
We estimate now the right hand side of this inequality. For the first term we have, using the estimates from Theorem \ref{characterization-decay-linear-part} 

\begin{align}
\label{eqn:estimate-linear-part}
 \int _0 ^t g'(s) \left( \int _{B(s)} \left( 1 + \a^2 |\xi|^2 \right) e^{2 t \mathcal{M} (\xi)} |\widehat{u_0} (\xi)|^2  \, d \xi \right) \, ds  \leq  C  \int _0 ^t g'(s) \Vert \overline{u} (s) \Vert ^2 _{H^1 _{\a}  (\R^3)}  \, ds \nonumber \\  \leq   C_1  \int _0 ^t g'(s) (C_2 + s) ^{- \left( \frac{3}{2} + r^{\ast} \right)}\, ds,
\end{align}
where $\overline{u}$ is the solution to the linear part of (\ref{eqn:navier-stokes-voigt}). For the second term, after integrating in polar coordinates in $B(t)$, we obtain

\begin{align}
\label{eqn:estimate-nonlinear-part}
 \int _0 ^t  g'(s) & \left( \int _{B(s)} \left(1  +  \a^2 |\xi|^2 \right)  |\xi|^2 \left( \int _0 ^t \Vert u(\tau) \Vert _{L^2 (\R^3)} ^2  \, d \tau \right)^2  \, d \xi \right) \, ds \nonumber \\ & \leq C \left( \int _0 ^t g'(s) \left( \rho^5 (s) + \rho^7 (s) \right)  \, ds \right) \left( \int_0 ^t \Vert u (\tau) \Vert _{L^2 (\R^3)} ^2 \, d \tau \right) ^2.
\end{align}
For a fixed $r^{\ast}$, we choose $g(t) = (C + t)^m$, for some $m > \max \{ \frac{1}{2},\frac{3}{2} + r^{\ast} \}$. Then $\rho(t) = C_1 (C_2 + t) ^{- \frac{1}{2}}$ and from (\ref{eqn:previous-to-main-inequality}), (\ref{eqn:estimate-linear-part}) and (\ref{eqn:estimate-nonlinear-part}) and the apriori estimate $\Vert u(t) \Vert _{L^2 (\R^3)} \leq C$, we obtain

\begin{eqnarray*}
\Vert u(t) \Vert ^2 _{H^1 _{\a} (\R^3)}  & \leq & C (1 + t) ^{-m} + C_1 (C_2 + t) ^{- \left( \frac{3}{2} + r^{\ast} \right)} + C_3 (C_4 + t) ^{- \frac{1}{2}} \\ & \leq & C_1 (C_2 + t) ^{- \min \{ \frac{3}{2} + r^{\ast}, \frac{1}{2}  \}}.
\end{eqnarray*}
We now use this first preliminary decay to bootstrap, trying to find sharper estimates for (\ref{eqn:estimate-nonlinear-part}). Suppose $ \min \{ \frac{3}{2} + r^{\ast}, \frac{1}{2}  \} = \frac{3}{2} + r^{\ast}$. Then, for $g(t) = (C + t)^m$ with $m > \max \{\frac{3}{2} + r^{\ast}, \frac{3}{2} \}$, we have that  $\rho(t) = C_1 (C_2 + t) ^{- \frac{1}{2}}$ and then

\begin{displaymath}
C \left( \int _0 ^t g'(s)  \left( \rho^5 (s) + \rho^7 (s) \right) ds \right) \left( \int_0 ^t \Vert u (\tau) \Vert _{L^2 (\R^3)} ^2 \, d \tau \right) ^2 \leq C_1 (C_2 + s) ^{- \left( \frac{7}{2} + 2 r^{\ast} \right)}.
\end{displaymath}
The decay in (\ref{eqn:estimate-linear-part}) is still the slower one, so there is no improvement. For $ \min \{ \frac{3}{2} + r^{\ast}, \frac{1}{2}  \} = \frac{1}{2}$, we obtain

\begin{displaymath}
 \left( \rho^5 (s) + \rho^7 (s) \right) \left( \int_0 ^s \Vert u (\tau) \Vert _{L^2 (\R^3)} ^2 \, d \tau \right) ^2 \leq C_1 (C_2 + s) ^{-  \frac{3}{2}}
\end{displaymath}
so going back to (\ref{eqn:estimate-nonlinear-part}) we have

\begin{eqnarray*}
\Vert u(t) \Vert _{H^1 _{\a} (\R^3)} ^2 & \leq & C (1 + t) ^{-m} + C_1 (C_2 + t) ^{- \left( \frac{3}{2} + r^{\ast} \right)} + C_3 (C_4 + t) ^{- \frac{3}{2}} \\ & \leq & C_1 (C_2 + t) ^{- \min \{ \frac{3}{2} + r^{\ast}, \frac{3}{2}  \}}.
\end{eqnarray*}
We boostrap once again. Let $ \min \{ \frac{3}{2} + r^{\ast}, \frac{3}{2}  \} = \frac{3}{2} + r^{\ast}$, and consider $r^{\ast} \neq - \frac{1}{2}$.  By the same computations as before, now with $m > \max \{\frac{3}{2} + r^{\ast}, \frac{5}{2} \}$, we show that decay is still given by  (\ref{eqn:estimate-linear-part}). If $r^{\ast} =  - \frac{1}{2}$, then the integral of the $L^2$ norm in (\ref{eqn:estimate-nonlinear-part}) leads to a logarithmic term which after choosing $g(t) = (C + t) ^m$, with $m > \frac{5}{2}$, leads to 

\begin{eqnarray*}
\Vert u(t) \Vert _{H^1 _{\a} (\R^3)} ^2 & \leq & C (1 + t) ^{-m} + C_1 (C_2 + t) ^{- 1} + C_3 \ln ^2 (1 + t) \, (C_3 + t) ^{- \frac{5}{2}} \\ & \leq & C (1 + t) ^{- 1}.
\end{eqnarray*}
If $ \min \{ \frac{3}{2} + r^{\ast}, \frac{3}{2}  \} = \frac{3}{2}$, then

\begin{displaymath}
\int _0 ^t \Vert u(\tau) \Vert _{L^2 (\R^3)} ^2 \, d \tau \leq C
\end{displaymath}
from which we obtain, after choosing $g(t) = (C + t) ^m$, with $m > \frac{5}{2}$

\begin{eqnarray*}
\Vert u(t) \Vert _{H^1 _{\a} (\R^3)} ^2 & \leq & C (1 + t) ^{-m} + C_1 (C_2 + t) ^{- \left( \frac{3}{2} + r^{\ast} \right)} + C_3 (C_4 + t) ^{- \frac{5}{2}} \\ & \leq & C (1 + t) ^{- \min \{\frac{3}{2} + r^{\ast},  \frac{5}{2} \}}.
\end{eqnarray*}
This finishes the proof of the Theorem. $\Box$

\smallskip

{\bf Proof of Theorem \ref{decay-nonlinear-part}} As before, we assume solutions are regular enough and we obtain estimates for them. The estimates ofr weak solutions are obtained as in Theorem \ref{main-theorem}. Take $w = u - \overline{u}$, where $\overline{u}$ is the solution to the linear part (\ref{eqn:linear-part}) with the same initial datum. Then $w$ solves the equation  

\begin{displaymath}
\partial_t (w - \a^2 \Delta w) + (u \cdot \nabla) u = \nu \Delta w - \nabla p, \quad w_0 (x) = 0.
\end{displaymath}
Multiplying by $w$ and integrating we obtain

\begin{eqnarray*}
\frac{1}{2} \frac{d}{dt} \left( \Vert w(t) \Vert _{L^2 (\R^3)} ^2 + \a ^2 \Vert \nabla w (t) \Vert _{L^2 (\R^3)} ^2  \right) + \nu \Vert \nabla w (t) \Vert _{L^2 (\R^3)} ^2 =  \int _{\R^3} \overline{u} (u \cdot \nabla) u \, dx.
\end{eqnarray*}
We now use interpolation for $L^3$ in terms of $L^2$ and $L^6$, the embedding $H^1 (\R^3) \subset L^6 (\R^3)$, the decay for $\overline{u}$ from Theorem \ref{characterization-decay-linear-part} and the decay for $u$ from Theorem \ref{main-theorem} to obtain

\begin{eqnarray}
\label{eqn:ineq-cross-term}
\int _{\R^3} \overline{u} (u \cdot \nabla) u \, dx & = & - \int _{\R^3} \nabla \overline{u} (u \otimes u) \, dx  \leq  \Vert \nabla \overline{u} (t) \Vert _{L^2 (\R^3)} \, \Vert u(t) \Vert _{L^3 (\R^3)} \Vert u (t) \Vert _{L^6 (\R^3)} \nonumber \\ & \leq &  \Vert \nabla \overline{u}(t) \Vert _{L^2 (\R^3)} \Vert u(t) \Vert ^{1/2} _{L^2 (\R^3)} \Vert u(t) \Vert ^{1/2} _{L^6 (\R^3)} \Vert  u (t) \Vert _{L^6 (\R^3)} \nonumber \\ & \leq & \Vert \overline{u} (t) \Vert _{H^1 _{\a} (\R^3)} \Vert u(t) \Vert ^{1/2} _{L^2 (\R^3)} \Vert u(t)  \Vert ^{3/2} _{H^1 _{\a} (\R^3)}  \nonumber \\ & \leq & \Vert  \overline{u} (t) \Vert _{H^1 _{\a} (\R^3)} \Vert u(t)  \Vert ^{2} _{H^1 _{\a} (\R^3)}  \nonumber \\ & \leq & C_1 (C_2 + t) ^{- \min \{ \frac{9}{4} + \frac{3}{2} r^{\ast}, \frac{13}{4} + \frac{1}{2} r^{\ast} \}}.
\end{eqnarray}

\br \label{remark-slow-decay-nsv} If we were working with the Navier-Stokes equations, we would have

\begin{displaymath}
\int _{\R^3} \overline{u} (u \cdot \nabla) u \, dx \leq  \Vert \nabla \overline{u} (t) \Vert _{L^{\infty} (\R^3)} \Vert u(t) \Vert ^2 _{L^2 (\R^3)} \leq C (1 + t) ^{- \frac{5}{4}} \Vert u(t) \Vert ^2 _{L^2 (\R^3)},
\end{displaymath}
instead of (\ref{eqn:ineq-cross-term}), see page 590 in Niche and M.E. Schonbek \cite{Niche01042015}. This is so because the heat kernel strongly regularizes initial data, unlike the linear part of the Navier-Stokes-Voigt equation (see Remark \ref{facts-linear-part-nsv}) and this allows us to estimate a higher order norm on $\nabla \overline{u}$ and to have some uniform decay. 
\er

Taking $B(t)$ as in (\ref{eqn:ball}) and following the same arguments as in the proof of Theorem \ref{main-theorem}, we have

\begin{eqnarray}
\label{eqn:inequality-with-integral}
\frac{d}{dt} \left( g(t) \Vert w(t) \Vert ^2 _{H^1 _{\a} (\R^3)} \right) &\leq& C g'(t) \int _{B(t)} \left( 1 + \a^2 |\xi|^2 \right) |\widehat{w} (\xi, t)| ^2 \, d \xi \nonumber \\ & + &  C_1 g'(t) (C_2 + t) ^{- \min \{ \frac{9}{4} + \frac{3}{2} r^{\ast}, \frac{13}{4} + \frac{1}{2} r^{\ast} \}}.
\end{eqnarray}
As

\begin{displaymath}
\widehat{w} (\xi,t) = \int _0 ^t e^{(t - s) \mathcal{M} (\xi) } \mathcal{F} ((u \cdot \nabla) u - \nabla p) (\xi,t) \, ds
\end{displaymath}
we obtain

\begin{displaymath}
|\widehat{w} (\xi, t)|^2 \leq C |\xi| \int_0 ^t \Vert u(s) \Vert  _{L^2 (\R^3)} ^2 \, ds,
\end{displaymath}
where we used that 

\begin{displaymath}
|\mathcal{F} \left( (u \cdot \nabla) u - \nabla p \right) (\xi, t)| \leq C |\xi| \Vert u(t) \Vert ^2 _{L^2 (\R^3)},
\end{displaymath}
which we proved in Lemma \ref{estimate-modulus-uhat}. Hence

\be
\label{eqn:estimate-ball}
\int _{B(t)} \left( 1 + \a^2 |\xi|^2 \right) |\widehat{w} (\xi, t)| ^2 \, d \xi \leq C \left(\rho^5 (t) + \rho ^7 (t)  \right) \left( \int_0 ^t \Vert u(s) \Vert  _{L^2 (\R^3)} ^2 \, ds \right) ^2.
\ee
We will choose $g$ so that $\rho(t) = C_1 (C_2+ t) ^{- \frac{1}{2}}$. Assuming we have done so, from (\ref{eqn:inequality-with-integral}) and (\ref{eqn:estimate-ball}) we obtain 

\begin{eqnarray*}
\label{eqn:inequality-with-g}
\frac{d}{dt} \left( g(t) \Vert w(t) \Vert ^2 _{H^1 _{\a} (\R^3) }\right) &\leq& C_1 g'(t)  (C_2+ t) ^{- \frac{5}{2}} \left( \int_0 ^t \Vert u(s) \Vert  _{L^2 (\R^3)} ^2 \, ds \right) ^2 \nonumber \\ & + & C_1 g'(t) (C_2 + t) ^{- \min \{ \frac{9}{4} + \frac{3}{2} r^{\ast}, \frac{13}{4} + \frac{1}{2} r^{\ast} \}}.
\end{eqnarray*}
If $r^{\ast} \neq - \frac{1}{2}$, then from the decay in Theorem \ref{main-theorem} we obtain

\begin{displaymath}
\left( \int_0 ^t \Vert u(s) \Vert  _{L^2 (\R^3)} ^2 \, ds \right)^2 \leq C_1 (C_2 + t) ^{- \min \{1 + 2 r^{\ast}, 3 \}}.
\end{displaymath}
which leads to 

\begin{displaymath}
\frac{d}{dt} \left( g(t) \Vert w(t) \Vert ^2 _{H^1 _{\a} (\R^3) }\right) \leq  C_1 g'(t) (C_2 + t) ^{- \min \{ \frac{9}{4} + \frac{3}{2} r^{\ast}, \frac{13}{4} + \frac{1}{2} r^{\ast}, \frac{11}{2}  \}}.
\end{displaymath}
Taking $g(t) = C_1 (C_2 + t) ^m$, with $m > \max \{ \frac{9}{4} + \frac{3}{2} r^{\ast}, \frac{13}{4} + \frac{1}{2} r^{\ast}, \frac{11}{2}  \}$ we obtain

\begin{displaymath}
\Vert w(t) \Vert ^2 _{H^1 _{\a} (\R^3)} \leq C_1 (C_2 + t) ^{- \min \{ \frac{9}{4} + \frac{3}{2} r^{\ast}, \frac{13}{4} + \frac{1}{2} r^{\ast}, \frac{11}{2}  \}}.
\end{displaymath}
If $r^{\ast} = - \frac{1}{2}$ we obtain a logarithmic term from the integral of the $L^2$ norm of u, which does not change the estimate obtained in the previous line. This finishes the proof. $\Box$

\bibliographystyle{plain}
\bibliography{NSV-Niche}

\end{document}